# Parameter Identifiability and Redundancy: Theoretical Considerations

Mark P. Little[1]*, Wolfgang F. Heidenreich[2], Guangquan Li[1]

1 Department of Epidemiology and Public Health, Imperial College Faculty of Medicine, St Mary's Campus, London, United Kingdom, 2 Institut für Strahlenschutz, Helmholtz Zentrum München, German Research Center for Environmental Health, Ingolstädter Landstrasse, Neuherberg, Germany

## Abstract

**Background:** Models for complex biological systems may involve a large number of parameters. It may well be that some of these parameters cannot be derived from observed data via regression techniques. Such parameters are said to be unidentifiable, the remaining parameters being identifiable. Closely related to this idea is that of redundancy, that a set of parameters can be expressed in terms of some smaller set. Before data is analysed it is critical to determine which model parameters are identifiable or redundant to avoid ill-defined and poorly convergent regression.

**Methodology/Principal Findings:** In this paper we outline general considerations on parameter identifiability, and introduce the notion of weak local identifiability and gradient weak local identifiability. These are based on local properties of the likelihood, in particular the rank of the Hessian matrix. We relate these to the notions of parameter identifiability and redundancy previously introduced by Rothenberg (*Econometrica* **39** (1971) 577–591) and Catchpole and Morgan (*Biometrika* **84** (1997) 187–196). Within the widely used exponential family, parameter irredundancy, local identifiability, gradient weak local identifiability and weak local identifiability are shown to be largely equivalent. We consider applications to a recently developed class of cancer models of Little and Wright (*Math Biosciences* **183** (2003) 111–134) and Little *et al.* (*J Theoret Biol* **254** (2008) 229–238) that generalize a large number of other recently used quasi-biological cancer models.

**Conclusions/Significance:** We have shown that the previously developed concepts of parameter local identifiability and redundancy are closely related to the apparently weaker properties of weak local identifiability and gradient weak local identifiability—within the widely used exponential family these concepts largely coincide.

**Funding:** This work was funded partially by the European Commission under contracts FI6R-CT-2003-508842 (RISC-RAD) and FP6-036465 (NOTE). The funders had no role in study design, data collection and analysis, decision to publish, or preparation of the manuscript.

**Competing Interests:** The authors have declared that no competing interests exist.

* E-mail: mark.little@imperial.ac.uk

## Introduction

Models for complex biological systems may involve a large number of parameters. It may well be that some of these parameters cannot be derived from observed data via regression techniques. Such parameters are said to be unidentifiable or non-identifiable, the remaining parameters being identifiable. Closely related to this idea is that of redundancy, that a set of parameters can be expressed in terms of some smaller set. Before data is analysed it is critical to determine which model parameters are identifiable or redundant to avoid ill-defined and poorly convergent regression.

Identifiability in stochastic models has been considered previously in various contexts. Rothenberg [1] and Silvey [2] (pp. 50, 81) defined a set of parameters for a model to be identifiable if no two sets of parameter values yield the same distribution of the data. Catchpole and Morgan [3] considered identifiability and parameter redundancy and the relations between them in a general class of (exponential family) models. Rothenberg [1], Jacquez and Perry [4] and Catchpole and Morgan [3] also defined a notion of local identifiability, which we shall extend in the Analysis Section. [There is also a large literature on identifiability in deterministic (rather than stochastic) models, for example the papers of Audoly *et al.* [5], and Bellu [6], which we shall not consider further.] Catchpole *et al.* [7] and Gimenez *et al.* [8] outlined use of computer algebra techniques to determine numbers of identifiable parameters in the exponential family. Viallefont *et al.* [9] considered parameter identifiability issues in a general setting, and outlined a method based on considering the rank of the Hessian for determining identifiable parameters; however, some of their claimed results are incorrect (as we outline briefly later). Gimenez *et al.* [8] used Hessian-based techniques, as well as a number of purely numerical techniques, for determining the number of identifiable parameters. Further general observations on parameter identifiability and its relation to properties of sufficient statistics are given by Picci [10], and a more recent review of the literature is given by Paulino and de Bragança Pereira [11].

In this paper we outline some general considerations on parameter identifiability. We shall demonstrate that the concepts of parameter local identifiability and redundancy are closely related to apparently weaker properties of weak local identifiability and gradient weak local identifiability that we introduce in the Analysis Section. These latter properties relate to the uniqueness of likelihood maxima and likelihood turning points within the vicinity of sets of parameter values, and are shown to be based on local properties of the likelihood, in particular the rank of the Hessian matrix. Within the widely-used exponential family we demonstrate that these concepts (local identifiability, redundancy, weak local identifiability, gradient weak local identifiability) largely coincide. We briefly consider applications of all these ideas to a recently developed general class of carcinogenesis models [12,13,14], presenting results that generalize those of Heidenreich [15] and Heidenreich *et al.* [16] in the context of the two-mutation cancer model [17]. These are outlined in the later parts of the Analysis and the Discussion, and in more detail in a companion paper [12].

## Analysis

### General Considerations on Parameter Identifiability

As outlined in the Introduction, a general criterion for parameter identifiability has been set out by Jacquez and Perry [4]. They proposed a simple linearization of the problem, in the context of models with normal error. They defined a notion of <u>local identifiability</u>, which is that in a local region of the parameter space, there is a unique $\theta_0$ that fits some specified body of data, $(x_i, y_i)_{i=1}^n$, i.e. for which the model predicted mean $h(x|\theta)$ is such that the residual sum of squares:

$$S = \sum_{l=1}^{n} [y_l - h(x_l|\theta)]^2 \quad (1)$$

has a unique minimum. We present here a straightforward generalization of this to other error structures. If the model prediction $h(x) = h(x|\theta)$ for the observed data $y$ is a function of some vector parameters $\theta = (\theta_j)_{j=1}^p$ then in general it can be assumed, under the general equivalence of likelihood maximization and iteratively reweighted least squares for generalized linear models [18](chapter 2) that one is trying to minimize:

$$S = \sum_{l=1}^{n} \frac{1}{v_l} \left[ y_l - h(x_l|\theta_0) - \sum_{j=1}^{p} \frac{\partial h(x_l|\theta)}{\partial \theta_j}\bigg|_{\theta=\theta_0} \cdot \Delta\theta_j \right]^2 \quad (2)$$

where $y_l (1 \leq l \leq n)$ $(n \geq p)$ is the observed measurement (e.g., the numbers of observed cases in the case of binomial or Poisson models) at point $l$ and the $v_l$ $(1 \leq l \leq n)$ are the current estimates of variance at each point. This has a unique minimum in the perturbing $\Delta\theta = (\Delta\theta_j)_{j=1}^p$ $(\theta = \theta_0 + \Delta\theta)$ given by $H^T D H \Delta\theta = H^T D \delta$, where $(\delta_l)_{l=1}^n = (y_l - h(x_l|\theta_0))_{l=1}^n$, $(H_{lj})_{l=1,j=1}^{n,p} = \left(\frac{\partial h(x_l|\theta)}{\partial \theta_j}\bigg|_{\theta=\theta_0}\right)_{l=1,j=1}^{n,p}$, $D = diag[1/v_1, 1/v_2, ..., 1/v_n]$, whenever $H^T D H$ has full rank $(=p)$.

More generally, suppose that the likelihood associated with observation $x_l$ is $l(x_l|\theta)$ and let $L(x_l|\theta) = \ln[l(x_l|\theta)]$. Then generalizing the least squares criterion (1) we now extend the definition of local identifiability to mean that there is at most one maximum of:

$$L = L(x|\theta) = \sum_{l=1}^{n} L(x_l|\theta) \quad (3)$$

in the neighborhood of any given $\theta \in \Omega \subset R^p$. More formally:

**Definitions 1.** A set of parameters $(\theta_i)_{i=1}^p$ is <u>identifiable</u> if for any $\theta \in \Omega$ there are no $\delta \in \Omega \setminus \{\theta\}$ for which $L(x|\delta) = L(x|\theta)$ ($x$ almost everywhere ($a.e.$)). A set of parameters $(\theta_i)_{i=1}^p$ is <u>locally identifiable</u> if there exists a neighborhood $N \in \aleph_\theta$ such that for no $\delta \in N \setminus \{\theta\}$ is $L(x|\delta) = L(x|\theta)$ ($x$ a.e.). A set of parameters $(\theta_i)_{i=1}^p$ is <u>weakly locally identifiable</u> if there exists a neighborhood $N \in \aleph_\theta$ and data $x = (x_1, ..., x_n) \in \Sigma^n$ such that the log-likelihood $L = L(x|\theta) = \sum_{l=1}^{n} L(x_l|\theta)$ is maximized by at most one set of $\hat{\theta} \in N$. If $L = L(x|\theta)$ is $C^1$ as a function of $\theta \in \Omega$ a set of parameters $(\theta_i)_{i=1}^p \in \text{int}(\Omega)$ is <u>gradient weakly locally identifiable</u> if there exists a neighborhood $N \in \aleph_\theta$ and data $x = (x_1, ..., x_n) \in \Sigma^n$ such that $\left(\frac{\partial L(x|\hat{\theta})}{\partial \hat{\theta}_i}\right)_{i=1}^p = 0$ (i.e., $\hat{\theta}$ is a <u>turning point</u> of $L(x|\theta)$) for at most one set of $\hat{\theta} \in N$.

Our definitions of identifiability and local identifiability coincide with those of Rothenberg [1], Silvey [2](pp. 50, 81) and Catchpole and Morgan [3]. Rothenberg [1] proved that if the Fisher information matrix, $I = I(\theta)$, in a neighborhood of $\theta \in \text{int}(\Omega)$ is of constant rank and satisfies various other more minor regularity conditions, then $\theta \in \text{int}(\Omega)$ is locally identifiable if and only if $I(\theta)$ is non-singular. Clearly identifiability implies local identifiability, which in turn implies weak local identifiability. By the Mean Value Theorem [19](p. 107) gradient weak local identifiability implies weak local identifiability. Heuristically, (gradient) weak local identifiability happens when:

$$0 = \frac{\partial L}{\partial \theta_i} = \sum_{l=1}^{n} \frac{\partial L(x_l|\theta)}{\partial \theta_i} = \sum_{l=1}^{n} \left[ \frac{\partial L(x_l|\theta_0)}{\partial \theta_i} + \sum_{j=1}^{p} \frac{\partial^2 L(x_l|\theta_0)}{\partial \theta_i \partial \theta_j} \cdot \Delta\theta_j \right] + O(|\Delta\theta|^2), 1 \leq i \leq p \quad (4)$$

and in general this system of $p$ equations has a unique solution in $\Delta\theta = (\Delta\theta_j)_{j=1}^p$ in the neighborhood of $\theta_0$ (assumed $\in \text{int}(\Omega)$) whenever $\left(\sum_{l=1}^{n} \frac{\partial^2 L(x_l|\theta_0)}{\partial \theta_i \partial \theta_j}\right)_{i,j=1}^p$ has full rank $(= p)$. This turns out to be (nearly) the case, and will be proved later (Corollary 2). More rigorously, we have the following result.

**Theorem 1.** Suppose that the log-likelihood $L(x|\theta)$ is $C^2$ as a function of the parameter vector $\theta \in \Omega \subset R^p$, for all $x = (x_1, ..., x_n) \in \Sigma^n$.

(i) Suppose that for some $x$ and $\theta \in \text{int}(\Omega)$ it is the case that $rk\left[\left(\frac{\partial^2 L(x|\theta)}{\partial \theta_i \partial \theta_j}\right)_{i,j=1}^p\right] = p$. Then turning points of the likelihood in the neighborhood of $\theta$ are isolated, i.e., there is an open neighborhood $N \in \aleph_\theta \subset \Omega$ for which there is at most one $\hat{\theta} \in N$ that satisfies $\left(\frac{\partial L(x|\theta)}{\partial \theta_i}\right)_{i=1}^p\bigg|_{\theta=\hat{\theta}} = 0$.

(ii) Suppose that for some $x$ and $\theta \in \text{int}(\Omega)$ it is the case that $rk\left[\left(\frac{\partial^2 L(x|\theta)}{\partial \theta_i \partial \theta_j}\right)_{i,j=1}^p\right] = p$ then local maxima of the likelihood in the neighborhood of $\theta$ are isolated, i.e., there is an open neighborhood $N \in \aleph_\theta \subset \Omega$ for which there is at most one $\hat{\theta} \in N$ that is a local maximum of $L(x|\theta)$.

(iii) Suppose that for some $x$ and all $\theta \in \text{int}(\Omega)$ it is the case that $rk\left[\left(\frac{\partial^2 L(x|\theta)}{\partial \theta_i \partial \theta_j}\right)_{i,j=1}^p\right] = r < p$ then all local maxima of the likelihood in $\text{int}(\Omega)$ are not isolated, as indeed are all $\theta \in \text{int}(\Omega)$ for which $\left(\frac{\partial L(x|\theta)}{\partial \theta_i}\right)_{i=1}^p = 0$.

We prove this result in Text S1 Section A. As an immediate consequence we have the following result.

**Corollary 1.** For a given $x=(x_1,...,x_n) \in \Sigma^n$, a sufficient condition for the likelihood (3) to have at most one maximum and one turning point in the neighborhood of a given $\theta = (\theta_1,...,\theta_p) \in \text{int}(\Omega)$ is that $rk\left[\left(\frac{\partial^2 L(x|\theta)}{\partial \theta_i \partial \theta_j}\right)_{i,j=1}^p\right] = p$. In particular, if this condition is satisfied $\theta$ is gradient weakly locally identifiable (and therefore weakly locally identifiable). ($\Omega \subset R^p$ is the parameter space.)

That this condition is not necessary is seen by consideration of the likelihood $l(x|\theta) = C \cdot \exp\left(-\sum_{i=1}^p [x_i - \theta_i]^4\right)$, where $C$ is chosen so that this has unit mass. Then $\frac{\partial^2 L(x|\theta)}{\partial \theta_i \partial \theta_j} = -12 \cdot [x_i - \theta_i]^2 \cdot \delta_{ij}$ which has rank 0 at $\theta = x$ and a unique maximum there. In particular, this shows that the result claimed by Viallefont et al. [9](proposition 2, p. 322) is incorrect.

**Definitions 2.** A subset of parameters $(\theta_{\pi(i)})_{i=1}^k$ (for some permutation $\pi: \{1,2,...,p\} \to \{1,2,...,p\}$) is <u>weakly maximal</u> (respectively <u>weakly gradient maximal</u>) if for any <u>permissible</u> fixed $(\theta_{\pi(i)})_{i=k+1}^p$ (such that $\Omega_{k,\pi}^{(\theta_{\pi(i)})_{i=k+1}^p} = \{(\theta_{\pi(i)})_{i=1}^k: (\theta_1,...,\theta_k,\theta_{k+1},...,\theta_p) \in \Omega\} \neq \varnothing$) $(\theta_{\pi(i)})_{i=1}^k$ is weakly locally identifiable (respectively gradient weakly locally identifiable) at that point, but that this is not the case for any larger number of parameters. A subset of parameters $(\theta_{\pi(i)})_{i=1}^k$ is <u>strongly maximal</u> (respectively <u>strongly gradient maximal</u>) if for any permissible fixed $(\theta_{\pi(i)})_{i=k+1}^p$ and any open $U \subset \Omega_{k,\pi}^{(\theta_{\pi(i)})_{i=k+1}^p}$, $(\theta_{\pi(i)})_{i=1}^k$ restricted to the set $U$ is weakly maximal (respectively weakly gradient maximal), i.e., all $(\theta'_{\pi(i)})_{i=1}^k \in U$ are weakly maximal (respectively weakly gradient maximal).

From this it easily follows that a strongly (gradient) maximal set of parameters $(\theta_{\pi(i)})_{i=1}^k$ is *a fortiori* weakly (gradient) maximal at all points $(\theta'_{\pi(i)})_{i=1}^k \in \Omega_{k,\pi}^{(\theta_{\pi(i)})_{i=k+1}^p}$ for any permissible $(\theta_{\pi(i)})_{i=k+1}^p$. Assume now that $k$ of the $p$ $\theta_i$ are a weakly maximal set of parameters. So for some permutation $\pi: \{1,2,...,p\} \to \{1,2,...,p\}$ and for any permissible fixed $(\theta_{\pi(i)})_{i=k+1}^p$ and any $(\theta_{\pi(i)})_{i=1}^k \in \Omega_{k,\pi}^{(\theta_{\pi(i)})_{i=k+1}^p} \subset R^k$ there is an open neighborhood $N \in \aleph_{(\theta_{\pi(i)})_{i=1}^k} \subset \Omega_{k,\pi}^{(\theta_{\pi(i)})_{i=k+1}^p}$ and some data $x=(x_1,...,x_n) \in \Sigma^n$ for which $L_{(\theta_{\pi(i)})_{i=k+1}^p}\left(x|(\theta_{\pi(i)})_{i=1}^k\right)$ is maximized by at most one set of $(\hat\theta_{\pi(i)})_{i=1}^k \in N$, but that this is not the case for any larger number of parameters. Assume that $r = \max\left\{rk\left[\left(\frac{\partial^2 L_{(\theta_{\pi(i)})_{i=k+1}^p}(x|(\theta_{\pi(i)})_{i=1}^k)}{\partial \theta_{\pi(i)} \partial \theta_{\pi(j)}}\right)_{i,j=1}^k\right]: (\theta_{\pi(i)})_{i=1}^k \in N\right\} < k$. If $L$ is $C^2$ as a function of $\theta$ then it follows easily that $\Omega_{k,r} = \left\{(\theta_{\pi(i)})_{i=1}^k \in N: rk\left[\left(\frac{\partial^2 L_{(\theta_{\pi(i)})_{i=k+1}^p}(x|(\theta_{\pi(i)})_{i=1}^k)}{\partial \theta_{\pi(i)} \partial \theta_{\pi(j)}}\right)_{i,j=1}^k\right] = r\right\}$ must be an open non-empty subset of $N$. By Theorem 1 (iii) any $\hat\theta \in \Omega_{k,r}$ which maximizes $L_{(\theta_{\pi(i)})_{i=k+1}^p}$ in $\Omega_{k,r}$ cannot be isolated, a contradiction (unless there are no maximizing $\hat\theta \in \Omega_{k,r}$). Therefore, either there are no maximizing $\hat\theta \in \Omega_{k,r}$ or for at least one $\hat\theta \in N$ $rk\left[\left(\frac{\partial^2 L_{(\theta_{\pi(i)})_{i=k+1}^p}(x|(\theta_{\pi(i)})_{i=1}^k)}{\partial \theta_{\pi(i)} \partial \theta_{\pi(j)}}\right)_{i,j=1}^k\bigg|_{(\theta_{\pi(i)})_{i=1}^k = \hat\theta}\right] = k$. This implies

that $rk\left[\left(\frac{\partial^2 L(x|\theta)}{\partial \theta_i \partial \theta_j}\right)_{i,j=1}^p\bigg|_{\theta = \hat\theta'}\right] \geq k$, where $\hat\theta' = (\hat\theta) \cup (\theta_{\pi(i)})_{i=k+1}^p$ in the obvious sense.

Assume now that the $(\theta_{\pi(i)})_{i=1}^k$ are strongly maximal. Suppose that for some $\theta_1 = (\theta_{1i})_{i=1}^p \in \Omega$ and some $x=(x_1,...,x_n) \in \Sigma^n$ it is the case that $rk\left[\left(\frac{\partial^2 L(x|\theta)}{\partial \theta_i \partial \theta_j}\right)_{i,j=1}^p\bigg|_{\theta = \theta_1}\right] > k$. Because $\left(\frac{\partial^2 L(x|\theta)}{\partial \theta_i \partial \theta_j}\right)_{i,j=1}^p\bigg|_{\theta = \theta_1}$ is symmetric, there is a permutation $\pi': \{1,...,p\} \to \{1,...,p\}$ for which $rk\left[\left(\frac{\partial^2 L(x|\theta)}{\partial \theta_{\pi'(i)} \partial \theta_{\pi'(i)}}\right)_{i,j=1}^{k+1}\bigg|_{\theta = \theta_1}\right] = k+1$ [20](p. 79). If $L$ is $C^2$ as a function of $\theta$ this will be the case in some open neighborhood $N' \in \aleph_{(\theta_{1\pi'(i)})_{i=1}^{k+1}} \subset R^{k+1}$. By Theorem 1 (ii) this implies that the parameters $(\theta_{\pi'(i)})_{i=1}^{k+1}$ have at most one maximum in $N'$, so that $(\theta_{\pi(i)})_{i=1}^k$ is not a strongly maximal set of parameters in $N'$. With small changes everything above also goes through with "weakly gradient maximal" substituted for "weakly maximal" and "strongly gradient maximal" substituted for "strongly maximal". Therefore we have proved the following result.

**Theorem 2.** Let $L(x|\theta)$ be $C^2$ as a function of $\theta \in \Omega \subset R^p$ for all $x \in \Sigma^n$.

(i) If there is a weakly maximal (respectively weakly gradient maximal) subset of $k$ parameters, $(\theta_{\pi(1)},\theta_{\pi(2)},...,\theta_{\pi(k)})$ (for some permutation $\pi: \{1,2,...,p\} \to \{1,2,...,p\}$), and for fixed $(\theta_{\pi(i)})_{i=k+1}^p$ and some $x=(x_1,...,x_n) \in \Sigma^n$ $L_{(\theta_{\pi(i)})_{i=k+1}^p}(x|(\theta_{\pi(i)})_{i=1}^k)$ has a maximum (respectively turning point) on the set of $\theta$ where $rk\left[\left(\frac{\partial^2 L_{(\theta_{\pi(i)})_{i=k+1}^p}(x|(\theta_{\pi(i)})_{i=1}^k)}{\partial \theta_{\pi(i)} \partial \theta_{\pi(j)}}\right)_{i,j=1}^k\right]$ is maximal then $\max\left\{rk\left[\left(\frac{\partial^2 L_{(\theta_{\pi(i)})_{i=k+1}^p}(x|(\theta_{\pi(i)})_{i=1}^k)}{\partial \theta_{\pi(i)} \partial \theta_{\pi(j)}}\right)_{i,j=1}^k\right] : (\theta_{\pi(i)})_{i=1}^k \in \Omega_{k,\pi}^{(\theta_{\pi(i)})_{i=k+1}^p}\right\} = k$ and $\max\left\{rk\left[\left(\frac{\partial^2 L(x|\theta)}{\partial \theta_i \partial \theta_j}\right)_{i,j=1}^p\right] : \theta \in \Omega\right\} \geq k$.

(ii) If there is a strongly maximal (respectively strongly gradient maximal) subset of $k$ parameters, $(\theta_{\pi(1)},\theta_{\pi(2)},...,\theta_{\pi(k)})$ (for some permutation $\pi: \{1,2,...,p\} \to \{1,2,...,p\}$) then $rk\left[\left(\frac{\partial^2 L(x|\theta)}{\partial \theta_{\pi(i)} \partial \theta_{\pi(j)}}\right)_{i,j=1}^k\right] \leq k \; \forall \theta \in \Omega$.

All further results in this Section assume that the model is a member of the exponential family, so that if the observed data $x=(x_l)_{l=1}^n \in \Sigma^n$ then the log-likelihood is given by $L(x|\theta) = \sum_{l=1}^n \left[\frac{x_l \varsigma_l - b(\varsigma_l)}{a(\phi)} + c(x_l,\phi)\right]$ for some functions $a(\phi), b(\varsigma), c(x,\phi)$. We assume that the natural parameters $\varsigma_l = \varsigma_l[(\theta_i)_{i=1}^p, z_l]$ are functions of the model parameters $(\theta_i)_{i=1}^p$ and some auxiliary data $z_l$, but that the scaling parameter $\phi$ is not. Let $\mu_l = b'(\varsigma_l) = E[x_l]$, so that $\mu_l = b'(\varsigma_l[(\theta_i)_{i=1}^p, z_l])$. In all that follows we shall assume that the function $b(\varsigma)$ is $C^2$. The following definition was introduced by Catchpole and Morgan [3].

**Definition 3.** With the above notation, a set of parameters $(\theta_i)_{i=1}^p \in \Omega$ is <u>parameter redundant</u> for an exponential family model if $\mu_l = b'(\varsigma_l[(\rho_i)_{i=1}^q, z_l])$ can be expressed in terms of some strictly smaller parameter vector $(\rho_i)_{i=1}^q$ $(q<p)$. Otherwise, the set of parameters $(\theta_i)_{i=1}^p$ is <u>parameter irredundant</u> or <u>full rank</u>.

Catchpole and Morgan [3] proved (their Theorem 1) that a set of parameters is parameter redundant if and only if $rk\left[\left(\frac{\partial \mu_l}{\partial \theta_i}\right)_{l=1,\,i=1}^{n\ \ p}\right]<p$. They defined full rank models to be <u>essentially full rank</u> if $rk\left[\left(\frac{\partial \mu_l}{\partial \theta_i}\right)_{l=1,\,i=1}^{n\ \ p}\right]=p$ for every $(\theta_i)_{i=1}^p \in \Omega$; if $rk\left[\left(\frac{\partial \mu_l}{\partial \theta_i}\right)_{l=1,\,i=1}^{n\ \ p}\right]=p$ only for some $(\theta_i)_{i=1}^p \in \Omega$ then the parameter set is <u>conditionally full rank</u>. They also showed (their Theorem 3) that if $I=I(\theta)$ is the Fisher information matrix then $rk\left[\left(\frac{\partial \mu_l}{\partial \theta_i}\right)_{l=1,\,i=1}^{n\ \ p}\right]=rk[I(\theta)]$, and that parameter redundancy implies lack of local identifiability; indeed their proof of Theorems 2 and 4 showed that there is also lack of weak local identifiability (respectively gradient weak local identifiability) for all $(\theta_i')_{i=1}^p \in \Omega$ which for some $x=(x_l)_{l=1}^n \in \sum^n$ are local maxima (respectively turning points) of the likelihood.

Assume that $\theta=(\theta_i)_{i=1}^p$ are an essentially full rank set of parameters for the model. From the above result for every $\theta=(\theta_i)_{i=1}^p \in \Omega$ $rk\left[\left(\frac{\partial \mu_l}{\partial \theta_i}\right)_{l=1,\,i=1}^{n\ \ p}\right]=rk(I(\theta))=p$. Therefore, since $E\left[\frac{\partial^2 L(x|\theta)}{\partial \theta_i \partial \theta_j}\right]=-E\left[\frac{\partial L(x|\theta)}{\partial \theta_i}\frac{\partial L(x|\theta)}{\partial \theta_j}\right]=-I(\theta)$ is of full rank and so negative definite, so by the strong law of large numbers we can choose $x=(x_l)_{l=1}^n \in \sum^n$ so that the same is true of $\frac{\partial^2 L(x|\theta)}{\partial \theta_i \partial \theta_j}=\sum_{l=1}^n\left\{\left[\frac{x_l-b'(\varsigma_l)}{a(\phi)}\right]\frac{\partial^2 \varsigma_l}{\partial \theta_i \partial \theta_j}-\frac{b''(\varsigma_l)}{a(\phi)}\frac{\partial \varsigma_l}{\partial \theta_i}\frac{\partial \varsigma_l}{\partial \theta_j}\right\}$. This implies that on some $N \in \aleph_\theta \subset R^p$ $\frac{\partial^2 L(x|\theta)}{\partial \theta_i \partial \theta_j}=\sum_{l=1}^n\left\{\left[\frac{x_l-b'(\varsigma_l)}{a(\phi)}\right]\frac{\partial^2 \varsigma_l}{\partial \theta_i \partial \theta_j}-\frac{b''(\varsigma_l)}{a(\phi)}\frac{\partial \varsigma_l}{\partial \theta_i}\frac{\partial \varsigma_l}{\partial \theta_j}\right\}$ is of full rank, and therefore by Corollary 1 $\theta=(\theta_i)_{i=1}^p$ is (gradient) weakly locally identifiable. Furthermore, the above argument shows that if $\theta=(\theta_i)_{i=1}^p$ are a conditionally full rank set of parameters then on the (open) set $\Omega_p=\left\{\theta=(\theta_i)_{i=1}^p \in \Omega: rk\left[\left(\frac{\partial \mu_l}{\partial \theta_i}\right)_{l=1,\,i=1}^{n\ \ p}\right]=p\right\}$ $\theta=(\theta_i)_{i=1}^p$ is gradient weakly locally identifiable. We have therefore proved:

**Theorem 3.** Let $L(x|\theta)$ belong to the exponential family and be $C^2$ as a function of $\theta \in \Omega \subset R^p$ for all $x \in \sum^n$.

(i) If the parameter set $\theta=(\theta_i)_{i=1}^p$ is parameter redundant then it is not locally identifiable, and is not weakly locally identifiable (respectively gradient weakly locally identifiable) for all $(\theta_i')_{i=1}^p \in \Omega$ which for some $x=(x_l)_{l=1}^n \in \sum^n$ are local maxima (respectively turning points) of the likelihood.

(ii) If the parameter set $\theta=(\theta_i)_{i=1}^p$ is of essentially full rank then for some $x=(x_l)_{l=1}^n \in \sum^n$ $\frac{\partial^2 L(x|\theta)}{\partial \theta_i \partial \theta_j}$ is of full rank and therefore $\theta=(\theta_i)_{i=1}^p$ is gradient weakly locally identifiable (and so weakly locally identifiable) for all $\theta=(\theta_i)_{i=1}^p \in \Omega$.

(iii) If the parameter set $\theta=(\theta_i)_{i=1}^p$ is of conditionally full rank then it is gradient weakly locally identifiable on the open set $\Omega_p=\left\{\theta=(\theta_i)_{i=1}^p \in \Omega: rk\left[\left(\frac{\partial \mu_l}{\partial \theta_i}\right)_{l=1,\,i=1}^{n\ \ p}\right]=p\right\}$.

Remarks: It should be noted that part (i) of this generalizes part (i) of Theorem 4 of Catchpole and Morgan [3], who proved that if a model is parameter redundant then it is not locally identifiable. However, some components of part (ii) (that being essentially full rank implies gradient weak local identifiability) is weaker than the other result, proved in part (ii) of Theorem 4 of Catchpole and Morgan [3], namely that if a model is of essentially full rank it is locally identifiable. As noted by Catchpole and Morgan [3] (pp. 193–4), there are exponential-family models that are conditionally full rank, but not locally identifiable, so part (iii) is about as strong a result as can be hoped for.

From Theorem 3 we deduce the following.

**Corollary 2.** Let $L(x|\theta)$ belong to the exponential family and be $C^2$ as a function of $\theta \in \Omega \subset R^p$ for all $x \in \sum^n$. Then

(i) If for some subset of parameters $(\theta_{\pi(i)})_{i=1}^k$ and some $x=(x_1,...,x_n) \in \sum^n$ it is the case that $rk\left[\left(\frac{\partial^2 L(x|\theta)}{\partial \theta_{\pi(i)} \partial \theta_{\pi(j)}}\right)_{i,j=1}^k\right]=k$ then this subset is gradient weakly locally identifiable at this point.

(ii) If a subset of parameters $(\theta_{\pi(i)})_{i=1}^k$ is weakly locally identifiable and for some $x \in \sum^n$ this point is a local maximum of the likelihood then it is parameter irredundant, i.e., of full rank, so $rk[I(\theta)]=k$, so that for some $x' \in \sum^{n'}$ $rk\left[\left(\frac{\partial^2 L(x'|\theta)}{\partial \theta_{\pi(i)} \partial \theta_{\pi(j)}}\right)_{i,j=1}^k\right]=k$. In particular, if this holds for all $\theta \in \Omega$ then parameter irredundancy, local identifiability, gradient weak local identifiability and weak local identifiability are all equivalent.

**Proof.** This is an immediate consequence of the remarks after Definition 1, Corollary 1, Theorem 3 (i) and Theorems 1 and 3 of Catchpole and Morgan [3]. **QED**.

Remarks: (i) By the remarks preceding Theorem 3 the conditions of part (i) (that for some $x=(x_1,...,x_n) \in \sum^n$ it is the case that $rk\left[\left(\frac{\partial^2 L(x|\theta)}{\partial \theta_i \partial \theta_j}\right)_{i,j=1}^k\right]=k$) are automatically satisfied if $\theta=(\theta_i)_{i=1}^k$ are an essentially full rank set of parameters for the model.

(ii) Assume the model is constructed from a stochastic cancer model embedded in the exponential family, in the sense outlined in Text S1 Section B, so that the natural parameters $\varsigma_l = \varsigma_l[(\theta_i)_{i=1}^p, z_l]$ are functions of the model parameters $(\theta_i)_{i=1}^p$ and some auxiliary data $(z_l)_{l=1}^n$, and the means are given by $\mu_l = b'(\varsigma_l[(\theta_i)_{i=1}^p, z_l]) = z_l \cdot h[(\theta_i)_{i=1}^p, y_l]$, where $h[(\theta_i)_{i=1}^p, y_l]$ is the cancer hazard function. In this case, as shown in Text S1 Section B,
$$\frac{\partial^2 L(x|\theta)}{\partial \theta_i \partial \theta_j}=\sum_{l=1}^n\left[\begin{array}{l}\frac{[x_l-b'(\varsigma_l)]z_l}{a(\phi)b''(\varsigma_l)}\frac{\partial^2 h(\theta,y_l)}{\partial \theta_i \partial \theta_j}\\-\frac{z_l^2}{a(\phi)}\frac{\partial h(\theta,y_l)}{\partial \theta_i}\frac{\partial h(\theta,y_l)}{\partial \theta_j}\left\{\frac{[b''(\varsigma_l)]^2+b'''(\varsigma_l)[x_l-b'(\varsigma_l)]}{[b''(\varsigma_l)]^3}\right\}\end{array}\right].$$
The second term inside the summation $\left(-\frac{z_l^2}{a(\phi)}\frac{\partial h(\theta,y_l)}{\partial \theta_i}\frac{\partial h(\theta,y_l)}{\partial \theta_j}\left\{\frac{[b''(\varsigma_l)]^2+b'''((\varsigma_l)[x_l-b'(\varsigma_l)]}{[b''(\varsigma_l)]^3}\right\}\right)_{i,j=1}^p$ is a rank 1 matrix and can be made small in relation to the first term, e.g., by making $z_l$ small. Therefore finding data $(x,y,z)=(x_1,...,x_n,y_1,...,y_n,z_1,...,z_n) \in \sum^n$ for which $rk\left[\left(\frac{\partial^2 L(x|\theta)}{\partial \theta_{\pi(i)} \partial \theta_{\pi(j)}}\right)_{i,j=1}^k\right]=k$ is equivalent to finding data for

which $rk\left[\left(\frac{\partial^2 h(\theta,y_l)}{\partial\theta_{\pi(i)}\partial\theta_{\pi(j)}}\right)_{i,j=1}^k\right]=k$, or by the result of Dickson [20](p. 79) for which $rk\left[\left(\frac{\partial^2 h(\theta,y_l)}{\partial\theta_i\partial\theta_j}\right)_{i,j=1}^p\right]=k$.

## Hessian vs Fisher Information Matrix as a Method of Determining Redundancy and Identifiability in Generalised Linear Models

We, as with Catchpole and Morgan [3], emphasise use of the Hessian of the likelihood rather than the Fisher information matrix considered by Rothenberg [1]. In the context of GLMs, we have $L(x|\theta) = \sum_{l=1}^n \left[\frac{x_l\varsigma_l - b(\varsigma_l)}{a(\phi)} + c(x_l,\phi)\right]$ and $g(\mu_i) = g(b'(\varsigma_i)) = \sum_{j=1}^p A_{ij}\theta_j$ for some link function $g()$ and fixed matrix $A$. We define $D_{ij} = \frac{\partial\mu_j}{\partial\theta_i} = \frac{1}{g'((\mu_j)A_{ji}} = (A^T G^{-1})_{ij}$ where $G = diag(g'(\mu_1), g'(\mu_2), ..., g'(\mu_n))$. Theorem 1 of Catchpole and Morgan [3] states that a model is parameter irredundant if and only if $rk[D]=p$. The score vector is given by $U_i = \frac{\partial L(x|\theta)}{\partial\theta_i} = \sum_{l=1}^n \frac{[x_l - \mu_l]}{a(\phi)} \frac{\partial\varsigma_l}{\partial\theta_i} = \sum_{l=1}^n \frac{[x_l - \mu_l]}{b''(\varsigma_l)a(\phi)} \frac{\partial\mu_l}{\partial\theta_i} = \frac{1}{a(\phi)}(D\Delta(x-\mu))_i$ where $\Delta = diag\left(\frac{1}{b''(\varsigma_1)}, \frac{1}{b''(\varsigma_2)}, ..., \frac{1}{b''(\varsigma_n)}\right)$. The Fisher information is therefore given by $I(\theta) = E[UU^T] = \frac{1}{a(\phi)^2} D\Delta V\Delta D^T$ where $V = \left(E\left[[x_i-\mu_i][x_j-\mu_j]\right]\right)_{i,j}$ is the data variance. Theorem 1 of Rothenberg [1] states that a model is locally identifiable if and only if $rk[I(\theta)]=p$. As above (Corollary 2 (ii)), heuristically parameter irredundancy, local identifiability, gradient weak local identifiability and weak local identifiability are all equivalent and occur whenever $rk(D\Delta V\Delta D^T) = rk(D) = p$. Clearly evaluating the rank of $D$ is generally much easier than that of $D\Delta V\Delta D^T$. Catchpole and Morgan [3] demonstrate use of Hessian-based methods to estimate parameter redundancy in a class of capture-recapture models.

However, for certain applications, both the Fisher information and the Hessian must be employed, as we now outline. Assume that the model is constructed from a stochastic cancer model embedded in an exponential family model in the sense outlined in Text S1 Section B. The key to showing that such an embedded model has no more than $N$ irredundant parameters is to construct (as is done in Little et al. [12]) some scalar functions $G_1(.), G_2(.), ..., G_N(.)$ such that the cancer hazard function $h(\theta)$ can be written as $h(G_1(\theta), G_2(\theta), ..., G_N(\theta))$. Since the cancer model is embedded in a member of the exponential family (in the sense outlined in Text S1 Section B) the same will be true of the total log-likelihood $L(x|\theta) = L(x|G_1(\theta), G_2(\theta), ..., G_N(\theta))$. By means of the Chain Rule we obtain $\frac{\partial^2 L(x|\theta)}{\partial\theta_i\partial\theta_j} = \sum_{l,k=1}^N \frac{\partial^2 L(x|G_1,...,G_N)}{\partial G_l\partial G_k}\frac{\partial G_l}{\partial\theta_i}\frac{\partial G_k}{\partial\theta_j} + \sum_{l=1}^N \frac{\partial L(x|G_1,...,G_N)}{\partial G_l}\frac{\partial^2 G_l}{\partial\theta_i\partial\theta_j}$, so that the Fisher information matrix is given by:

$$I(\theta) = -E_\theta\left[\frac{\partial^2 L(x|\theta)}{\partial\theta_i\partial\theta_j}\right] = -E\left[\sum_{l,k=1}^N \frac{\partial^2 L(x|G_1,...,G_N)}{\partial G_l\partial G_k}\frac{\partial G_l}{\partial\theta_i}\frac{\partial G_k}{\partial\theta_j}\right] \quad (5)$$
$$= -\sum_{l,k=1}^N \frac{\partial G_l}{\partial\theta_i} E\left[\frac{\partial^2 L(x|G_1,...,G_N)}{\partial G_l\partial G_k}\right]\frac{\partial G_k}{\partial\theta_j}$$

which therefore has rank at most $N$. Therefore by Corollary 2 there can be at most $N$ irredundant parameters, or indeed (gradient) weak locally identifiable parameters. [A similar argument shows that if one were to reparameterise (via some invertible $C^2$ mapping $\theta = f(\omega)$) then the embedded log-likelihood $L(x|f^{-1}(\theta)) = L(x|\omega)$ associated with $h(f^{-1}(\theta)) = h(\omega)$ must also have Fisher information matrix of rank at most $N$.] By remark (ii) after Corollary 2, to show that a subset of cardinality $N$ of the parameters $(\theta_i)_{i=1}^p$ is (gradient) weak locally identifiable parameters, requires that one show that $\left[\frac{\partial^2 h(\theta,y_l)}{\partial\theta_i\partial\theta_j}\right]_{i,j=1}^p$ has rank at least $N$ for some $(\theta,y_l)$. This is the approach adopted in the paper of Little et al. [12].

## Discussion

In this paper we have introduced the notions of weak local identifiability and gradient weak local identifiability, which we have related to the notions of parameter identifiability and redundancy previously introduced by Rothenberg [1] and Catchpole and Morgan [3]. In particular we have shown that within the exponential family models parameter irredundancy, local identifiability, gradient weak local identifiability and weak local identifiability are largely equivalent.

The slight novelty of our approach is that the notions of weak local identifiability and gradient weak local identifiability that we introduce are related much more to the Hessian of the likelihood rather than the Fisher information matrix that was considered by Rothenberg [1]. However, in practice, the two approaches are very similar; Catchpole and Morgan [3] used the Hessian of the likelihood, as do we, because of its greater analytic tractability. The use of this approach is motivated by the application, namely to determine identifiable parameter combinations in a large class of stochastic cancer models, as we outline at the end of the Analysis Section. In certain applications the Fisher information may be best for estimating the upper bound to the number of irredundant parameters, but the Hessian may be best for estimating the lower bound of this quantity.

In the companion paper of Little et al. [12] we consider the problem of parameter identifiability in a particular class of stochastic cancer models, those of Little and Wright [13] and Little et al. [14]. These models generalize a large number of other quasi-biological cancer models, in particular those of Armitage and Doll [21], the two-mutation model [17], the generalized multistage model of Little [22], and a recently developed cancer model of Nowak et al. [23] that incorporates genomic instability. These and other cancer models are generally embedded in an exponential family model in the sense outlined in Text S1 Section B, in particular when cohort data are analysed using Poisson regression models, e.g., as in Little et al. [13,14,24]. As we show at the end of the Analysis Section, proving (gradient) weak local identifiability of a subset of cardinality $k$ of the parameters $(\theta_i)_{i=1}^p$ can be done by showing that for this subset of parameters $rk\left[\left(\frac{\partial^2 h(\theta,y)}{\partial\theta_i\partial\theta_j}\right)_{i,j=1}^p\right]=k$ where $h$ is the cancer hazard function. Little et al. [12] demonstrate (by exhibiting a particular parameterization) that there is redundancy in the parameterization for this model: the number of theoretically estimable parameters in the models of Little and Wright [13] and Little et al. [14] is at most two less than the number that are theoretically available, demonstrating (by Corollary 2) that there can be no more than this number of irredundant parameters. Two numerical examples suggest that this bound is sharp – we show that the rank of the Hessian, $rk\left[\left(\frac{\partial^2 h(\theta,y)}{\partial \theta_i \partial \theta_j}\right)_{i,j=1}^p\right]$, is two less than the row dimension of this matrix. This result generalizes previously derived results of Heidenreich and others [15,16] relating to the two-mutation model.

## Supporting Information

### Text S1

Found at: doi:10.1371/journal.pone.0008915.s001 (0.33 MB DOC)

## Acknowledgments

The authors are very grateful for the comments of Professor Byron Morgan on an advanced draft of the paper, also for the detailed and helpful remarks of a referee.

## Author Contributions

Conceived and designed the experiments: MPL WFH. Performed the experiments: MPL GL. Analyzed the data: MPL GL. Wrote the paper: MPL WFH GL.

1 **Supplementary material A.**

2 **Proof of Theorem 1**

3 In this Section we outline a proof of Theorem 1 in the main text. To prove this result we need

4 the following lemma of Rudin [19](p.229).

5 **Lemma A1.** Suppose $m, n, r$ are non-negative integers s.t. $m \geq r$, $n \geq r$ and $F$ is a $C^1$

6 function $E \subset R^n \to R^m$ where $E$ is an open set. Suppose that $rk(F'(x)) = r$ $\forall x \in E$. Fix

7 $a \in E$ and put $A = F'(a)$, and let $Y_1 = A(R^n)$ and let $P: R^m \to R^m$ be a linear projection

8 operator ($P^2 = P$) s.t. $Y_1 = P(R^m)$ and $Y_2 = null(P)$. Then $\exists U, V \subset R^n$, open sets and a

9 bijective $C^1$ function $H: V \to U$ whose inverse is also $C^1$ and s.t. $F(H(x)) = Ax + \varphi(Ax)$,

10 $\forall x \in V$ where $\varphi: AV \subset Y_1 \to Y_2$ is a $C^1$ function.

11 We now restate Theorem 1 here.

12 **Theorem A2.** Suppose that the log-likelihood $L(x|\theta)$ is $C^2$ as a function of the parameter

13 vector $\theta \in \Omega \subset R^p$, and for all $x = (x_1, ..., x_n) \in \Sigma^n$.

14 (i) Suppose that for some $x$ and $\theta \in \text{int}(\Omega)$ it is the case that $rk\left[\left(\frac{\partial^2 L(x|\theta)}{\partial \theta_i \partial \theta_j}\right)_{i,j=1}^p\right] = p$.

15 Then turning points of the likelihood in the neighborhood of $\theta$ are isolated, i.e., there is

16 an open neighborhood $N \in \aleph_\theta \subset \Omega$ for which there is at most one $\hat{\theta} \in N$ that satisfies

17 $\left(\frac{\partial L(x|\theta)}{\partial \theta_i}\right)_{i=1}^p \bigg|_{\theta=\hat{\theta}} = 0$.

18 (ii) Suppose that for some $x$ and $\theta \in \text{int}(\Omega)$ it is the case that $rk\left[\left(\frac{\partial^2 L(x|\theta)}{\partial \theta_i \partial \theta_j}\right)_{i,j=1}^p\right] = p$

19 then local maxima of the likelihood in the neighborhood of $\theta$ are isolated, i.e., there is an

20  open neighborhood $N \in \aleph_\theta \subset \Omega$ for which there is at most one $\hat{\theta} \in N$ that is a local

21  maximum of $L(x|\theta)$.

22  (iii) Suppose that for some $x$ and all $\theta \in \text{int}(\Omega)$ it is the case that

23  $$rk\left[\left(\frac{\partial^2 L(x|\theta)}{\partial \theta_i \partial \theta_j}\right)_{i,j=1}^p\right] = r < p$$ then all local maxima of the likelihood in $\text{int}(\Omega)$ are not

24  isolated, as indeed are all $\theta \in \text{int}(\Omega)$ for which $\left(\frac{\partial L(x|\theta)}{\partial \theta_i}\right)_{i=1}^p = 0$.

25  **Proof:**

26  (i) Let $F: \Omega \subset R^p \to R^p$ be defined by $F(\theta_1, \theta_2, ..., \theta_p) = \left(\frac{\partial L(x|\theta)}{\partial \theta_1}, \frac{\partial L(x|\theta)}{\partial \theta_2}, ..., \frac{\partial L(x|\theta)}{\partial \theta_p}\right)$.

27  Since $L$ is $C^2$, $F$ is $C^1$ on $\text{int}(\Omega) \subset R^p$. By assumption $\frac{\partial F(\theta)_i}{\partial \theta_j} = \left(\frac{\partial^2 L(x|\theta)}{\partial \theta_j \partial \theta_i}\right)_{i,j=1}^p$ is of full

28  rank at $\theta$. By the inverse function theorem [19](pp.221-223) there are open $N, M \subset R^p$ such

29  that $\theta \in N$ and a $C^1$ bijective function $G: M \to N$ such that $G(F(\hat{\theta})) = \hat{\theta}$ for all $\hat{\theta} \in N$. In

30  particular there can be at most a single $\hat{\theta} \in N$ for which $F(\hat{\theta}) = 0$. **QED.**

31  (ii) By (i) there is an open neighborhood $N \in \aleph_\theta \subset \Omega$ for which if $\hat{\theta} \in N$ is such that

32  $\left(\frac{\partial L(x|\theta)}{\partial \theta_i}\right)_{i=1}^p\bigg|_{\theta=\hat{\theta}} = 0$ then for $\theta' \neq \hat{\theta} \in N$ $\left(\frac{\partial L(x|\theta)}{\partial \theta_i}\right)_{i=1}^p\bigg|_{\theta=\theta'} \neq 0$. Suppose now that $\hat{\theta} \in N$ is

33  a local maximum of $L(x|\theta)$. Any member of this neighborhood other than $\hat{\theta}$ cannot be a

34  turning point, and so by the Mean Value Theorem (Rudin 1976, p.107) cannot be a local

35  maximum. **QED.**

36     (iii) Let $F: \Omega \subset R^p \to R^p$ be defined by

37     $F(\theta_1, \theta_2, ..., \theta_p) = \left( \frac{\partial L(x|\theta)}{\partial \theta_1}, \frac{\partial L(x|\theta)}{\partial \theta_2}, ..., \frac{\partial L(x|\theta)}{\partial \theta_p} \right)$. Since $L$ is

38     $C^2$, $F$ is $C^1$ on $\text{int}(\Omega) \subset R^p$. By assumption $rk(F) = rk\left[ \left( \frac{\partial^2 L(x|\theta)}{\partial \theta_i \partial \theta_j} \right)_{i,j=1}^p \right] = r$ for all

39     $\theta \in \text{int}(\Omega) \subset R^p$. Suppose that $\theta_0 \in \text{int}(\Omega)$ is a local maximum of $L$. Let

40     $A = \left. \frac{\partial F}{\partial \theta} \right|_{\theta=\theta_0} : R^p \to R^p$ ($A \in L(R^p, R^p)$), and choose some arbitrary projection

41     $P \in L(R^p, R^p)$ s.t. $P(R^p) = Y_1 = A(R^p)$, and let $Y_2 = null(P)$. By Lemma A1 there are open

42     sets $U, V \subset R^p$ with $\theta_0 \in U \subset \text{int}(\Omega)$ and a bijective $C^1$ mapping with $C^1$ inverse

43     $H: V \to U$ s.t. $F(y) = AH^{-1}y + \varphi(AH^{-1}y)$, $\forall y \in V$ where $\varphi: AV \subset Y_1 \to Y_2$ is a $C^1$

44     function.

45     Since $\theta_0 \in \text{int}(\Omega)$ is a local maximum of $L(x;\theta)$, by the Mean Value Theorem [19](p.107)

46     $F(\theta_0) = 0$. Now choose some non-trivial vector $\underset{\sim}{k} \in null(A)$ and define a function, as we can,

47     on some interval $\delta: (-\varepsilon, \varepsilon) \to R^p$ by $\delta(t) = H(H^{-1}(\theta_0) + t\underset{\sim}{k})$. Because $H: V \to U$ is bijective

48     and $\underset{\sim}{k}$ is non-trivial $\delta(t) = \delta(t') \Leftrightarrow t = t'$. Also, it is the case that:

49     $$\begin{aligned} F(\delta(t)) &= AH^{-1}(H[H^{-1}(\theta_0) + t\underset{\sim}{k}]) + \varphi(AH^{-1}(H[H^{-1}(\theta_0) + t\underset{\sim}{k}])) = \\ &A[H^{-1}(\theta_0) + t\underset{\sim}{k}] + \varphi(A[H^{-1}(\theta_0) + t\underset{\sim}{k}]) = A[H^{-1}(\theta_0)] + \varphi(A[H^{-1}(\theta_0)]) = F(\theta_0) = 0 \end{aligned} \quad (A1)$$

50     Define $G: (-\varepsilon, \varepsilon) \to R$ by $G(t) = L(\delta(t)) = L((\delta_1(t), \delta_2(t), ..., \delta_n(t)))$. By the chain rule

51     [19](p.215) $\frac{dG}{dt} = \sum_{i=1}^p \frac{\partial L(x|\theta)}{\partial \theta_i} \frac{d\delta_i}{dt} = 0$ $\forall t \in (-\varepsilon, \varepsilon)$. Finally, by the Mean Value Theorem

52     [19](p.107) $G$ must be constant; in particular $L(x|\delta(t)) = L(x|\delta(0)) = L(x|\theta_0)$ $\forall t \in (-\varepsilon, \varepsilon)$

53     and so all points $\delta(t)$ must also be local maxima of $L(x|\theta)$. Therefore $\theta_0$ is not an isolated

54   local maximum. Since all we used about $\theta_0 \in \text{int}(\Omega)$ was that

55   $F(\theta_0) = 0$ $F((\theta_{0i})_{i=1}^{p}) = \left( \dfrac{\partial L(x|\theta_0)}{\partial \theta_{01}}, \dfrac{\partial L(x|\theta_0)}{\partial \theta_{02}}, ..., \dfrac{\partial L(x|\theta_0)}{\partial \theta_{0p}} \right) = 0$, the above argument also

56   shows that turning points cannot be isolated: $F(\delta(t)) = 0$. **QED.**

57

## Supplementary material B.

## Specification of embedded exponential family model

In this Section we outline the specification of an embedding of a stochastic cancer model in a general class of statistical models, the so-called exponential family [18]. This is often done in fitting cancer models to epidemiological and biological data (e.g., see references [12, 13, 14, 24]). Recall that a model is a member of the <u>exponential family</u> if the observed data

$x = (x_l)_{l=1}^n \in \Sigma^n$ is such that the log-likelihood is given by $L(x|\theta) = \sum_{l=1}^n \left[ \frac{x_l \varsigma_l - b(\varsigma_l)}{a(\phi)} + c(x_l, \phi) \right]$

for some functions $a(\phi), b(\varsigma), c(x, \phi)$. We assume that the natural parameters $\varsigma_l = \varsigma_l[(\theta_i)_{i=1}^p, z_l]$ are functions of the model parameters $(\theta_i)_{i=1}^p$ and some auxiliary data $(z_l)_{l=1}^n$, and that $\mu_l = b'(\varsigma_l[(\theta_i)_{i=1}^p, z_l]) = z_l \cdot h[(\theta_i)_{i=1}^p, y_l]$. Here $h[(\theta_i)_{i=1}^p, y_l]$ is the cancer hazard function (for example, that of Little *et al.* [14], as also specified in the main text and in Text S1 Section B of Little *et al.* [12]), $(y_l)_{l=1}^n$ are some further auxiliary data, and we assume that the $(z_l)_{l=1}^n$ are all non-zero. [<u>Note:</u> this is not necessarily a generalized linear model (GLM).]

In this case it is seen that

$$\frac{\partial^2 L(x|\theta)}{\partial \theta_i \partial \theta_j} = \sum_{l=1}^n \left[ \frac{[x_l - b'(\varsigma_l)] z_l}{a(\phi) b''(\varsigma_l)} \frac{\partial^2 h(\theta, y_l)}{\partial \theta_i \partial \theta_j} - \frac{z_l^2}{a(\phi)} \frac{\partial h(\theta, y_l)}{\partial \theta_i} \frac{\partial h(\theta, y_l)}{\partial \theta_j} \left\{ \frac{[b''(\varsigma_l)]^2 + b'''(\varsigma_l)[x_l - b'(\varsigma_l)]}{[b''(\varsigma_l)]^3} \right\} \right] \quad (B1)$$

so that the Fisher information matrix is given by

$$I(\theta)_{ij} = -E_\theta \left[ \frac{\partial^2 L(x|\theta)}{\partial \theta_i \partial \theta_j} \right] = \frac{1}{a(\phi)} \sum_{l=1}^n \frac{z_l^2}{b''(\varsigma_l)} \frac{\partial h(\theta, y_l)}{\partial \theta_i} \frac{\partial h(\theta, y_l)}{\partial \theta_j} \quad (B2)$$